\documentclass[11pt]{amsart}
\usepackage{amsthm,amstext,amsmath,amscd,amssymb,latexsym}
\usepackage[matrix,arrow]{xy}

\date{}
 \newtheorem{thm}{Theorem}[section]

 \newtheorem{prop}[thm]{Proposition}
\theoremstyle{definition}
 
 \newtheorem{ex}[thm]{Example}
 \newtheorem{rem}[thm]{Remark}

\def\H      {\mathbb{H}}                % Quaternionen

\def\P      {\mathbb{P}}                % projektiver Raum

\def\EE  {\mathcal E}

\def\FF  {\mathcal F}
\def\GG  {\mathcal G}
\def\MM  {\mathcal M}

\def\OO  {\mathcal O}

\def\QQ  {\mathcal Q}

\def\H   {\mathrm  H}

\def\ker{\mathop{\mathrm{Ker}}\nolimits}
\def\im{\mathop{\mathrm{Im}}\nolimits}
\def\coker{\mathop{\mathrm{Coker}}\nolimits}

\def\th{\mathop{\mathrm{th}}\nolimits}

\def\Syz{\mathop{\mathrm{Syz}}\nolimits}

\numberwithin{equation}{section}

\begin{document}

%\begin{frontmatter}
\title[Kumar's correspondence]
{Implementation of Kumar's Correspondence }
\author{Hirotachi Abo}
\address{Colorado State University, Dept. of Mathematics, Fort
Collins, CO 80523, USA} 
\email{abo@math.colostate.edu}
\author{Chris Peterson}
\address{Colorado State University, Dept. of Mathematics, Fort
Collins, CO 80523, USA} 
\email{peterson@math.colostate.edu}
%\subjclass{Primary 14F05, 14Q99; Secondary 13C10} 
%\end{keyword}


%\end{frontmatter}

%\topmargin -1.5cm
%\textwidth 15.5cm
%\textheight 23cm
%\oddsidemargin 0.25cm
%\evensidemargin -1.0cm

\date{}
%-------------------------------------------------------------------
%\setcounter{section}{0}
%\renewcommand  {\sectionmark}[1]   {\markright{\thesection\ #1}}

 %\newtheorem{thm}{Theorem}[section]
 %\newtheorem{lem}[thm]{Lemma}
% \newtheorem{cor}[thm]{Corollary}
% \newtheorem{prop}[thm]{Proposition}
%\theoremstyle{definition}
% \newtheorem{df}[thm]{Definition}
 %\newtheorem{ex}[thm]{Example}
% \newtheorem{rem}[thm]{Remark}

%\begin{document}

\maketitle

\begin{abstract}
In 1997, N.M. Kumar published a paper which introduced a new tool of
use in the construction of algebraic vector bundles. Given a vector
bundle on projective $n$-space, a well known theorem of Quillen-Suslin guarantees
the existence of sections which generate the bundle on the
complement of a hyperplane in projective $n$-space. Kumar used this fact to give a
correspondence between vector bundles on projective $n$-space and vector bundles
on projective $(n-1)$-space satisfying certain conditions. He then applied this
correspondence to establish the existence of many, previously
unknown, rank two bundles on projective fourspace in positive characteristic. The
goal of the present paper is to give an explicit homological
description of Kumar's correspondence in a setting appropriate for
implementation in a computer algebra system.
\end{abstract}

\section{Introduction}
A fundamental problem in algebraic geometry is the study,
classification and construction of varieties, schemes and sheaves.
These problems are related in the sense that progress in one area
often leads to progress in each of the other areas. For instance,
given a sheaf with interesting or unusual properties, one can often
obtain correspondingly interesting varieties and schemes as
degeneracy loci of the sheaf. A main focus of the present paper is
an explicit homological description of a tool of use in the
construction of locally free sheaves on $\P^n$ over an algebraically
closed field, $K$, of arbitrary characteristic. With a slight abuse
of language, we will use the term Algebraic Vector Bundle for such a
sheaf. A vector bundle $\EE$ of rank $r$ on $\P^n$ is said to be of
low rank if $r<n$. The co-rank of a bundle is the difference $n-r$.
It appears that indecomposable low rank vector bundles on $\P^n$ are
exceedingly rare. In fact, the only known co-rank 2 vector bundles
in characteristic zero are the Horrocks-Mumford bundle on $\P^4$ and
the Horrocks bundle on $\P^5$ \cite{horrocks2,HM}. In characteristic
$p>2$ there are the additional co-rank 2 constructions of Kumar, and
Kumar et al \cite{kumar,KPR}. In characteristic $p=2$ there is a
single example of an indecomposable co-rank 3 bundle constructed by
Tango \cite{tango}. It is an open problem to construct other
examples or show that they do not exist. In particular, it is
unknown if there exist co-rank 2, indecomposable vector bundles on
$\P^n$ for any value of $n$ greater than 5.  An interesting class of
problems is concerned with establishing the existence or
non-existence of higher co-rank bundles on $\P^n$ with prescribed
properties. \vskip .1in
 The first constructions of higher co-rank algebraic vector bundles appeared
 in the 1970's in the papers of Horrocks-Mumford, Horrocks and
 Tango.
After Horrock's paper in 1978, no fundamentally new, higher co-rank
bundles were shown to exist for 20 years. In 1997, Kumar introduced
a completely novel construction method and demonstrated its power by
constructing several previously unknown co-rank 2 vector bundles in
positive characteristic \cite{kumar}. His method provided fuel for
the additional constructions found in \cite{KPR}. Kumar based his
construction on the solution, by Quillen and Suslin, of the
well-known Serre's conjecture on the existence of finitely
generated, non-free $K[x_0, \cdots , x_n]$-modules
\cite{serre,quillen,suslin}. For a given vector bundle on the
$n$-dimensional projective space $\P^n$, the theorem of Quillen and
Suslin guarantees us the existence of sections that generate the
vector bundle on the complement of a hyperplane in $\P^n$. The pair
of the vector bundle and these sections corresponds to a vector
bundle on the hyperplane. Kumar gave necessary and sufficient
conditions for a vector bundle on a hyperplane of $\P^n$ to be
obtained from a vector bundle on $\P^n$ in this way. His
correspondence between vector bundles on $\P^n$ and vector bundles
on a hyperplane (satisfying certain conditions) were used to
establish the existence of many, previously unknown, rank two vector
bundles on $\P^4$ in positive characteristic.

 The
purpose of the present paper is to give an explicit homological
description of Kumar's correspondence in a setting appropriate for
implementation in a computer algebra system.
\section{Preliminaries}

\subsection{Kumar's correspondence}

Let $K$ be a field. In 1955, J.P. Serre asked whether there exist
finitely generated $K[x_0, \cdots , x_n]$-modules which are not
free \cite{serre}. In 1976, Quillen and Suslin independently
proved that such modules do not exist, i.e. they showed that every
finitely generated projective $K[x_0, \cdots , x_n]$-module is
free (cf.~\cite{quillen}, \cite{suslin}). One can apply the
theorem of Quillen and Suslin to vector bundles on $\P^n$ as
follows. Let $h$ be a linear form in $K[x_0, \cdots , x_n]$.
 Let~$H$ be the hyperplane in~$\P^n$ determined by the zeros of
 $h$. Let $\EE$ be a vector bundle on $\P^n$ of rank $r$.
 By the theorem of Quillen and Suslin, $\EE$ restricted to
 the complement, $\P^n \setminus H$ of~$H$, is free.
 As a consequence, there exist $r$
sections~$s_1,\dots ,s_r \in \H^0 \left(\P^n \setminus H, \EE^{\vee}
|_{\P^n \setminus H}\right)$ that
 generate $\EE^{\vee}|_{\P^n \setminus H}$.
 It is known that for suitable integers~$l_i$,
 $1 \leq i \leq r$, the sections~$h^{l_i}s_i$
 extend to global sections~$\widetilde{s_i} \in \H^0(\P^n, \EE^{\vee}(l_i))$
 (cf.~\cite{hartshorne2}).
 Such sections define an injective morphism
 of sheaves~$\EE \rightarrow \bigoplus^r_{i = 1} \OO_{\P^n}(l_i)$,
 which is an injective bundle map outside
 the divisor defined by~$\widetilde{s_1} \wedge \cdots \wedge \widetilde{s_r}
 \in \H^0 \left(\P^n, (\wedge^r\EE^{\vee})(\sum^r_{i=1}l_i)\right)
 \cong \H^0 \left(\P^n, \OO_{\P^n}(\sum^r_{i=1}l_i - c_1(\EE))\right)$.
 By construction, this divisor is the $m^{th}$
 infinitesimal neighborhood~$H_m$ of~$H$,
 where~$m = \sum^r_{i=1}l_i - c_1(\EE)$.
 In other words,
 there is an exact sequence
\begin{equation} \label{eq:exact_sequence_2.1}
 0 \rightarrow \EE \rightarrow \bigoplus^r_{i = 1} \OO_{\P^n}(l_i) \rightarrow
 \FF \rightarrow 0,
\end{equation}
where $\FF$ is a coherent sheaf whose support is $H$. It is clear
that the coherent sheaf $\FF$ on $\P^n$ possesses an
$\OO_{H_m}$-module structure, and from~(\ref{eq:exact_sequence_2.1})
it follows that the homological dimension of $\FF$ is 1. Conversely,
if there exists a coherent sheaf $\FF$ on $\P^n$ which has an
$\OO_{H_m}$-structure, has homological dimension $1$ and which
allows a surjective morphism from a direct sum of $r$ line bundles
then there exists a rank $r$ vector bundle $\EE$ on $\P^n$ and an
exact sequence of type~(\ref{eq:exact_sequence_2.1}).

Let $\pi$ be the finite morphism $\pi: H_m \rightarrow H$ induced
by the projection~$\P^n \setminus P \rightarrow H$ from a point~$P
\in \P^n \setminus H$. Then $\pi_*$ induces an equivalence of
categories from the category of quasi-coherent $\OO_{H_m}$-modules
to the category of quasi-coherent $\OO_H$-modules having a $\pi_*
\OO_{H_m}$-module structure. This correspondence enables us to
translate statements about quasi-coherent $\OO_{H_m}$-modules into
statements about quasi-coherent $\OO_H$-modules.

\vspace{2mm} \noindent (1) Since $\pi_* \OO_{H_m} \simeq
\bigoplus_{i = 0}^{m-1}\OO_H(-i)$, a quasi-coherent $\OO_H$-module
$\QQ$ has a $\pi_* \OO_{H_m}$-module structure if and only if there
is a morphism $\phi : \QQ \rightarrow \QQ(1)$ whose $m^{th}$ power
is zero. Following Kumar, we call such a morphism a {\it nilpotent
endomorphism of} $\QQ$. From the theorem of Auslander and Buchsbaum
it follows that a quasi-coherent $\OO_{H_m}$-module has homological
dimension $1$ as a coherent sheaf on $\P^n$ if and only if the
corresponding quasi-coherent $\OO_H$-module has homological
dimension $0$, in other words, if the $\OO_H$-module is a vector
bundle.

\vspace{2mm} \noindent (2) Let $\MM$ be the direct image sheaf of
$\FF$ by $\pi$ and let $\phi$ be the corresponding nilpotent
endomorphism of $\MM$. Since $\pi$ is a finite morphism, there are
natural isomorphisms $\H^0(\P^n,\FF(-l_i)) \simeq \H^0(H,
\MM(-l_i))$ for all $1 \leq i \leq r$. We denote the restriction
of $\GG=\bigoplus^r_{i = 1} \OO_{\P^n}(l_i)$ to $H$ by $\GG_H$.
There is a surjective morphism from $\GG$ to $\FF$ if and only if
the restriction map from $\bigoplus_{i=0}^{m-1} \GG_H(-i)$ to
$\MM$ is surjective. The latter condition is equivalent to the
condition that there exists a map~$\psi:\GG_H\rightarrow \MM$ such
that $(\phi, \psi):\MM(-1) \oplus \GG_H\rightarrow\MM$ is
surjective.

\begin{thm}(Kumar)
There is a correspondence between {\rm (i)} and {\rm (ii)}:
\begin{itemize}
\item[{\rm (i)}] The set of pairs $(\EE, s)$,
where $\EE$ is a rank $r$ vector bundle on $\P^n$ and $s$ is a
morphism from $\EE$ to $\bigoplus^r_{i=1} \OO(l_i)$ with cokernel
$\FF$ satisfying:

a) $\FF$ is a coherent sheaf on the $m^{th}$ infinitesimal
neighborhood $H_m$ of a hyperplane $H$ for some positive integer
$m$.

b) The direct image sheaf of $\FF$ by the finite morphism $\pi : H_m
\rightarrow H$ is a vector bundle $\MM$ on $H$.

 \

\item[{\rm (ii)}] The set of triples $(\MM, \phi, \psi)$,
where $\MM$ is a vector bundle on $H$, $\phi:\MM\rightarrow
\MM(1)$ is a nilpotent endomorphism and $\psi:\bigoplus_{i=1}^{r}
\OO_H(l_i)\rightarrow\MM$ is a morphism such that $(\phi, \psi)
:\MM(-1)\oplus \bigoplus_{i=1}^{r} \OO_H(l_i)\rightarrow\MM$ is
surjective.
\end{itemize}
\label{th:kumar-corr}
\end{thm}
\begin{proof}
See~\cite{kumar} for a detailed proof.
\end{proof}
Our goal is to make explicit the procedure for computing the pair
$(\EE,s)$ corresponding to a given triple $(\MM,\phi,\psi)$ and
conversely, to make explicit the procedure for computing the triple
$(\MM,\phi,\psi)$ corresponding to a given pair $(\EE,s)$. Let $R$
be the homogeneous coordinate ring of $\P^{n-1}$ and $S$ the
homogeneous coordinate ring of $\P^n$. Suppose that there exists a
morphism $s$ from a rank $r$ vector bundle $\EE$ on $\P^n$ to
$\bigoplus^r_{i=1} \OO_{\P^n}(l_i)$ satisfying the condition in
Theorem~\ref{th:kumar-corr}. Then $s$ induces a homomorphism from
$\H^0_*(\P^n, \EE)$ to $\H^0_*(\P^n,\bigoplus^r_{i=1}
\OO_{\P^n}(l_i))=\bigoplus^r_{i=1} S(l_i)$. The sheafification of
the cokernel of $s$ is the sheaf $\FF$. From the cokernel of $s$ we
can compute the module $F=\H^0_* (\P^n,\FF)$. Consider the
$R$-module $_RF$ obtained from $F$ by restriction of scalars. Then
the sheaf associated to $_RF$ is $\MM$. So the key step in each
procedure is to compute the $R$-module $_RF$ from an $S$-module $F$
or an $S$-module $F$ from an $R$-module $M$ such that $_RF=M$. In
the following section we will discuss how to carry out these steps.

\subsection{Restriction of scalars}
Let $S$ be the polynomial ring $K[x_0, \dots ,x_n]$ and let $R$ be
the polynomial ring $K [x_0,\dots,x_{n-1}]$. For any graded
$S$-module $F$ we denote by $_RF$ the $R$-module obtained from $F$
by restriction of scalars. Let $Q$ be the quotient ring
$S/(x_n^m)$ for some integer $m$. Suppose that $F$ is finitely
generated and has a $Q$-module structure (i.e. $F$ is annihilated
by the ideal $(x_n^m)$). Then $_RF$ is also finitely generated and
has an $_RQ$-module structure. Indeed, the following proposition
immediately follows from the definition of restriction of scalars.
\begin{prop}
Let $F$ be a finitely generated graded $S$-module with minimal
generating set $\mathfrak{F}=\{f_i\}_{1 \leq i \leq s}$. Suppose
that $F$ has a $Q$-module structure. Then $\mathfrak{M} =
\{x_n^if_j \}_{{0 \atop 1}{\leq \atop \leq } {i \atop j}{\leq
\atop \leq}{m-1 \atop s}}$ is a generating set for $_RF$. Moreover
the $_RQ$-module structure of $_RF$ is determined by the
homomorphism $\phi : \mbox{$_RF$} \rightarrow \mbox{$(_RF)(1)$}$
defined by
\[
 x_n^if_j \mapsto
  \begin{cases}
    0 & i \geq m-1  \\
    x_n^{i+1}f_j & \mbox{otherwise}.
  \end{cases}
\]
\end{prop}

\begin{rem}
(1) The homomorphism~$\phi : \mbox{$_RF$} \rightarrow
\mbox{$(_RF)(1)$}$ corresponds to multiplication~$\cdot x_n : F
\rightarrow F(1)$, and clearly the {\it m}$^{\th}$ power of $\phi$
is zero. The homomorphism $\phi : \mbox{$_RF$} \rightarrow
\mbox{$(_RF)(1)$}$ obtained in this way will be called {\it the
standard nilpotent endomorphism} of $_RF$.

\noindent (2) The generating set $\mathfrak{M}$ of $_RF$ is not
always minimal. Eliminating redundant elements gives a minimal set
$\mathfrak{M}'=\{g_1, \dots, g_t\}$ of generators for $_RF$. Let
\[
M_0 \rightarrow \mbox{$_RF$} \rightarrow 0
\]
be the corresponding epimorphism, where $M_0$ is a free
$R$-module. Note that each $x_n g_i$ can be written as an
$R$-linear combination of $g_1, \dots ,g_t$:
\[
 x_n g_i = \sum_{j = 1}^t a_{ij}g_j,
\]
where $a_{ij} \in R$. So the matrix $(a_{ij})_{1 \leq i,j \leq t}$
defines a lifting $\phi_0 : M_0 \rightarrow M_0(1)$ of the
standard nilpotent endomorphism $\phi$ of $M$, since $\phi$ sends
$g_i$ to $x_n g_i$ for $1 \leq i \leq t$. We call the
lifting~$\phi_0$ of $\phi$ given in this way {\it the standard
lifting of $\phi$}. \label{th:rem-nilpotent}
\end{rem}

A homomorphism from a finitely generated $R$-module $M$ to $M(1)$ is
said to be a {\it nilpotent endomorphism} of $M$ if its $m^{th}$
power is zero for some positive integer $m$. The functor $_R \cdot$
induces an equivalence of categories from the category
$\mathfrak{S}_m$ of finitely generated $S$-modules having a
$Q=S/(x_n^m)$-module structure to the category $\mathfrak{R}$ of
finitely generated $R$-modules having an $_RQ$-module structure
(i.e. having a nilpotent endomorphism $\phi$ with $\phi^m=0$).
Indeed, for an $R$-module $M=<g_1,\dots,g_t>$, we can define a
finitely generated $S$-module $^S M$ by considering the set of all
$S$-linear combinations of the generators of $M$ (i.e. the set
$\{b_1g_1+ \cdots + b_tg_t \ | \ b_i \in S \}$). Its $Q$-module
structure is defined by
\begin{equation}
 \phi(g_i)=x_ng_i  \
\mbox{for each $i=1, \dots ,t$. } \label{eq:def-str}
\end{equation}
Obviously the functors $_R \cdot$ and $^S \cdot$ are inverse to
each other.

For each $i=1, \dots, t$, $x_ng_i$ can be written as an $R$-linear
combination of the $g_j$'s by (\ref{eq:def-str}), so we can define
the {\it standard lifting} for $\phi$ in the same way as in
Remark~\ref{th:rem-nilpotent}. The following proposition will show
us how to compute from $M$ the corresponding module $^S M$:

\begin{prop}
Let $M$ be an object of $\mathfrak{R}$ and let $\phi$ be a
nilpotent endomorphism of $M$ with $\phi^m=0$. Suppose that $M$
has a minimal free presentation of type
\begin{equation} M_1 \stackrel{\alpha}{\rightarrow} M_0 \rightarrow
M \rightarrow 0 \label{eq:min-free-pre}
\end{equation}
Then the corresponding $S$-module $F$ in $\mathfrak{S}_m$ has a
presentation
\[
 \begin{CD}
\left(M_1 \otimes_R S \right) \oplus \left(M_0(-1) \otimes_R S
\right) @>\left(\alpha, \hskip 1pt \phi_0(-1) - \cdot\hskip 1pt x_n
\right)>> M_0 \otimes_R S \rightarrow F \rightarrow 0
 \end{CD}
 \]

\noindent where $\phi_0 : M_0 \rightarrow M_0(1)$ is the standard
lifting of~$\phi$ and $\cdot x_n$
%:  M_0(-1) \otimes_R S \rightarrow M_0 \otimes_R S$
is multiplication by~$x_n$. \label{th:free_presentation}
\end{prop}
\begin{proof}
Let $\{g_1, \dots, g_t\}$ be a minimal set of generators for $M$.
Then $F=\{b_1g_1+\cdots+b_tg_t \ | \ b_i \in S \}$. Let $\phi_0
(-1)=(a_{ij})_{1\leq i,j \leq t}$ be the standard lifting of
$\phi(-1)$. Then it follows from (\ref{eq:def-str}) that
$\{g_1,\dots,g_t\}$ satisfies the relations
\begin{equation}
\sum_{j=1}^t a_{ij}g_j-x_ng_i=0 \label{eq:relation}
\end{equation}
for all $i=1, \dots,t$. So $(\alpha, \phi_0(-1)-\cdot x_n)$ forms
part of a presentation matrix of $F$. Suppose that there is a
relation on $\{g_1, \dots,g_t\}$:
\[
 c_1g_1+\cdots+c_tg_t=0,
\]
where $c_i \in S$ for each $i$. Without loss of generality, we may
assume that each term $c_ig_i$ can be rewritten in the form
$(c_i'x_n+c_i'')g_i$, where $c_i' \in S$ and $c_i'' \in R$. Let
$C=(c_1, c_2, \dots, c_t)$, $C'=(c_1', c_2', \dots, c_t')$,
$C''=(c_1'', c_2'', \dots, c_t'')$, $G=(g_1,g_2,\dots ,g_t)$ and
$A=(a_{ij})$. By using the relations given in (\ref{eq:relation}),
we get $c_1g_1+\cdots+c_tg_t=CG^T=C'AG^T+C''G^T$. Set
$b_j=\Sigma_{i=1}^tc_i'a_{ij}+c_j''$ and $B=[b_1,b_2,\dots ,b_j]$.
Then $CG^t=BG^T$. View $c_i,b_i$ as elements of $R[x_n]$. Let
$r=max\{deg(c_i)|\hskip 1pt 1\leq i\leq t\}$ and
$s=max\{deg(b_i)|\hskip 1pt 1\leq i\leq t\}$. The construction
guarantees that $s<r$. If we now repeat the same operation with
$b_1g_1+\cdots+b_tg_t$ then in a finite number of steps we can
decrease the maximum degree of the coefficients of the syzygy until
all of the coefficients have degree 0, i.e. the relation becomes an
$R$-linear combination of the $g_i$ which is equal to $0$:
\[
 d_1g_1+ \cdots +d_tg_t=0, \ d_i \in R \  \mbox{for each $i$}.
\]
Since we assumed that the presentation of $M$ given in
(\ref{eq:min-free-pre}) is minimal, $(d_1,\dots,d_t)^t$ can be
generated by column vectors of $\alpha$. Therefore,
$(\alpha,\phi_0(-1)-\cdot x_n)$ is a presentation matrix of $F$.
\end{proof}

\section{Algorithm}

In this section we will develop a procedure for computing a rank $r$
vector bundle on $\P^n$ from a given vector bundle on $\P^{n-1}$
satisfying the conditions in Theorem~\ref{th:kumar-corr}. The
procedure takes as input a triple $(\MM, \phi, \psi)$ and produces
as output the corresponding pair $(\EE,s)$. More specifically, the
procedure takes as input:
\begin{itemize}
\item  The finitely generated $R$-module $M=\ <g_1, \dots ,g_t>$
with minimal free presentation
\[ M_1 \stackrel{\alpha}{\rightarrow} M_0 \rightarrow M
\rightarrow 0\] whose associated sheaf, $\MM=\widetilde{M}$, is
locally free;
\item A nilpotent endomorphism $\phi$ of $M$ and its standard lifting $\phi_0$;
\item A homomorphism
$\psi=(\psi_1, \dots, \psi_r)$ from a free module $\bigoplus_{i=1}^r
R(l_i)$ to $M$ such that the corresponding sheaf morphism from
$\bigoplus_{i=1}^r \OO(l_i)$ to $\MM$  is a morphism such that
$(\phi, \psi) :\MM(-1)\oplus \bigoplus_{i=1}^{r}
\OO_H(l_i)\rightarrow\MM$ is surjective.
\end{itemize}
The procedure produces as output:
\begin{itemize}
\item The finitely generated $S$-module $E$
whose associated sheaf is a rank $r$ vector bundle;
\item A homomorphism $s : E \rightarrow \bigoplus_{i=1}^r S(l_i)$
such that the coherent sheaf associated to $\coker (s)$ coincides
with $^S M$.
\end{itemize}
To get the pair $(\EE,s)$ from the triple $(\MM,\phi,\psi)$, we take
the following steps:
\begin{itemize}
\item[(i)] Define a finitely generated $S$-module $F$ by
$\{a_1g_1+\cdots+a_tg_t \ | \ a_i \in R\}$. In practice, this
module will be given as the cokernel of the homomorphism $(\alpha,
\phi_0(-1)-\cdot x_n ):(M_1 \otimes_R S) \oplus (M_0(-1) \otimes_R
S)\rightarrow M_0 \otimes_R S$ (see
Proposition~\ref{th:free_presentation}).
\item[(ii)] Define the homomorphism from $\bigoplus_{i=1}^r S(l_i)$ to $F$
by $\psi=(\psi_1, \dots, \psi_r)$ and compute the syzygy module
$\Syz (\psi_1, \dots ,\psi_r)$ which represents the desired
homomorphism $s : E \rightarrow \bigoplus_{i=1}^r S(l_i)$. Note
that $\psi_i$ can be written as an $R$-linear combination of the
$g_j$'s for each $i=1, \dots, t$. So a simple way of computing
$\Syz (\psi_1, \dots ,\psi_r)$ is to determine the generating set
$\{g_1, \dots ,g_t\}$ of $F$ as a $Q=S/(x_n^m)$-module by using
the presentation matrix of $F$ given in (i). This enables us to
compute $\Syz (\psi_1, \dots ,\psi_r)$ as a $Q$-module. Indeed,
let $N$ be the extension of the module $\Syz(\psi_1, \dots
,\psi_r)$ to $S$. Then $\Syz (\psi_1, \dots ,\psi_r)$ will be
obtained as the quotient of $N$ by $x_n^m N$.
\end{itemize}
\begin{rem}
Let $(E,s)$ be the resulting pair. Then we want to check that
$\EE=\widetilde{E}$ is indeed a rank $r$ vector bundle on $\P^n$.
By construction, $\EE$ can be regarded as a subsheaf of
$\bigoplus_{i=1}^r \OO_{\P^n}(l_i)$:
\[
\begin{array}{ccrclcc}
\cdots & \rightarrow & \bigoplus_{j=1}^k \OO_{\P^n}(m_j) &
\xrightarrow{A} & \bigoplus_{i=1}^r \OO_{\P^n}(l_i) & \rightarrow
& \cdots \\
& & \searrow & & \nearrow & & \\
& & & \EE & & & \\
& & \nearrow & & \searrow & & \\
& & 0\ \ \ \ \  & & \ \ \ \ \ 0  & &
\end{array}
\]

%$$
%\xymatrix{ \cdots \ar[r] & \bigoplus_{j=1}^k \OO_{\P^n}(m_j)
%\ar[dr] \ar[rr]^{A} &&  \bigoplus_{i=1}^r \OO_{\P^n}(l_i)
%\ar[r] & \cdots ,\\
%&& \EE \ar[dr] \ar[ur]^{s} \\
%&0 \ar[ur] && 0 }
%$$
The entries of the $j^{th}$ column of $A$ define the scheme of
zeros $X_{s_j}=\{s_j=0\}$; the entries of the $i^{th}$ row of $A$
define the scheme of zeros $X_{\sigma_i}=\{\sigma_i=0\}$. Recall
that $s$ is an injective bundle map outside the divisor defined by
\[
 x_n^m=\sigma_1 \wedge \cdots \wedge \sigma_r
 \in \H^0(\P^n,(\wedge^r\EE^{\vee})(\sum^r_{i=1}l_i))
 \cong \H^0(\P^n, \OO_{\P^n}(m)),
\]
where $c_1$ is the first Chern class of $\EE$ and
$m=\sum_{i=1}^rl_i -c_1$. The $j^{th}$ column of $A$ represents
the section $t_j=s(s_j)$ of
$\bigoplus_{i=1}^r\OO_{\P^n}(l_i-m_j)$. So we have the relation of
the form
\[
t_{j_1} \wedge \cdots \wedge t_{j_r}=x_n^m \cdot (s_{j_1} \wedge
\cdots \wedge s_{j_r}),
\]
and hence we can prove that $\widetilde{E}$ is a vector bundle by
checking that the ideal quotient $(I:x_n^m)$ defines the empty set
in $\P^n$, where $I$ is the ideal generated by the maximal minors
of $A$. \label{th:locally-free}
\end{rem}

The following examples will show how the procedure works. The
procedure in the first example takes as input the twisted cotangent
bundle on $\P^2$ and returns as output a stable rank two vector
bundle on $\P^3$ with Chern classes $(c_1,c_2)=(0,1)$. This bundle
is the null correlation bundle on $\P^3$.
\begin{ex}
Let $R=K[x_0, x_1, x_2]$ and let $S=K[x_0,  x_1,x_2 ,x_3]$.
Consider the following Koszul complex:
\[ 0\rightarrow R(-1) \stackrel{\alpha_2}{\longrightarrow} 3R
\stackrel{\alpha_1}{\longrightarrow} 3R(1)
\stackrel{\alpha_0}{\longrightarrow} R(2) \] where
\[
\alpha_0=(x_0,x_1,x_2), \ \ \alpha_1= \left(
\begin{array}{ccc}
-x_1 & -x_2 & 0 \\
x_0 & 0 & -x_2 \\
0 & x_0 & x_1
\end{array}
\right) \ \ \mbox{and} \ \ \alpha_2= \left(
\begin{array}{c}
x_2 \\
-x_1 \\
x_0
\end{array}
\right).
\]
Let $M=\mathrm{Im} (\alpha_1)=\ <s_1,s_2,s_3>$. Then $\widetilde{M}$
is the twisted cotangent bundle $\Omega^1(2)$. The third row, $t_1$
of $\alpha_1$, induces a map from $\Omega^1 (2)$ to $\OO(1)$ such
that $t_1 \circ s_1=0$. So the composite of $s_1(1)$ and $t_1$
defines a nilpotent endomorphism $\phi$ of $M$, and hence
$\widetilde{M}$. In this case, the standard lifting of $\phi$ is
\[
 \phi_0=
\left(
\begin{array}{ccc}
0 & x_0 & x_1 \\
0 & 0 & 0 \\
0 & 0 & 0
\end{array}
\right): 3R \rightarrow 3R(1).
\]
This can be summarized in the following sequence of maps
\[
 \dots\rightarrow R\stackrel{{\tiny\left(\begin{array}{c} 1\\ 0\\ 0\end{array} \right)}}{\longrightarrow} 3R
 \stackrel{{\tiny\left(
\begin{array}{ccc}
-x_1 & -x_2 & 0 \\
x_0 & 0 & -x_2 \\
0 & x_0 & x_1
\end{array}
\right)}}{\longrightarrow} 3R(1) \stackrel{{\tiny\left(\begin{array}{ccc}
0 & 0 & 1\end{array} \right)}}{\longrightarrow}
 R(1)\stackrel{{\tiny\left(\begin{array}{c} 1\\ 0\\ 0\end{array} \right)}}{\longrightarrow} 3R(1)
 \rightarrow \dots
\]

The fact that $t_1 \circ s_1=0$ corresponds to \[
\left(\begin{array}{ccc} 0 & 0 & 1\end{array} \right)\left(
\begin{array}{ccc}
-x_1 & -x_2 & 0 \\
x_0 & 0 & -x_2 \\
0 & x_0 & x_1
\end{array}
\right)\left(\begin{array}{c} 1\\ 0\\ 0\end{array} \right)=0.\]
The map $\phi_0:3R\rightarrow 3R(1)$ corresponds to
 \[
 \phi_0=
\left(
\begin{array}{ccc}
0 & x_0 & x_1 \\
0 & 0 & 0 \\
0 & 0 & 0
\end{array}
\right)=\left(\begin{array}{c} 1\\ 0\\ 0\end{array}
\right)\left(\begin{array}{ccc} 0 & 0 & 1\end{array} \right)\left(
\begin{array}{ccc}
-x_1 & -x_2 & 0 \\
x_0 & 0 & -x_2 \\
0 & x_0 & x_1
\end{array}
\right).
\]

By Proposition~\ref{th:free_presentation}, the corresponding
$S$-module $F$ in $\mathfrak{S}_2$ has the following minimal
presentation:
\[
4S(-1) \stackrel{\beta_0}{\longrightarrow} 3S \rightarrow F
\rightarrow 0,
\]
where the first column of $\beta_0$ is the presentation matrix for
$M$ (i.e. $\alpha_2$) and the next three columns of $\beta_0$ are
just the columns of the matrix $\phi_0(-1)-x_3I$ where I is the
$3\times 3$ identity matrix. Thus,
\[
 \beta_0 =
\left(
\begin{array}{cccc}
x_2 & -x_3 & x_0 & x_1 \\
-x_1 & 0 & -x_3 & 0 \\
x_0 & 0 & 0 & -x_3
\end{array}
\right).
\]
The other generators $s_2$ and
$s_3$ of $M$ define a homomorphism $\psi : 2R \rightarrow M$,
whose lifting is given by the matrix
\[
\psi_0= \left(
\begin{array}{cc}
0 & 0 \\
1 & 0 \\
0 & 1 \\
\end{array}
\right): 2R\rightarrow 3R.
\]
This homomorphism together with the nilpotent endomorphism
$\phi(-1)$ of $M(-1)$ yields a homomorphism $(\phi(-1), \psi):
M(-1) \oplus 2R \rightarrow M$. The image $N$ is generated by the
columns of the matrix
\[
\left(
\begin{array}{ccccc}
0 & -x_0x_1 & -x_1^2 & -x_2 & 0 \\
0 & x_0^2 & x_0x_1 & 0 & -x_2 \\
0 & 0  & 0 &  x_0 & x_1
\end{array}
\right): 3R(-1)\oplus 2R \rightarrow 3R(1).
\]
The first three columns of the matrix come from $\alpha_1\phi_0(-1)$
(i.e. multiply $\alpha_1$ and $\phi_0$) and the next two columns
come from $\alpha_1\psi_0$ (i.e. multiply $\alpha_1$ and $\psi_0$).
The truncated modules $M_{\geq 1}$ and $N_{\geq 1}$ are isomorphic,
so the map of sheaves $(\phi(-1),\psi): \Omega^1(1) \oplus 2\OO
\rightarrow \Omega^1(2)$ is surjective. From
Theorem~\ref{th:kumar-corr} it follows that there exists a rank two
vector bundle $\EE$ on $\P^3$ with exact sequence

\begin{equation}
 0 \rightarrow \EE \rightarrow 2\OO \rightarrow \widetilde{F} \rightarrow 0.
\label{sq:rank-two-bundle}
\end{equation}

Let $\mathfrak{F}=\{f_1,f_2,f_3\}$ be the minimal generating set
of $F$, where for each $i$, $f_i$ corresponds to $s_i$. By
construction, the surjective map from $2\OO$ to $\widetilde{F}$ in
Sequence~(\ref{sq:rank-two-bundle}) is induced by $f_2$ and $f_3$.
Let $Q$ be the quotient ring $S/(x_3^2)$. Then $F$, as a
$Q$-module, is generated by
\[
f_1= \left(
\begin{array}{c}
0 \\
x_1x_3 \\
x_0x_3
\end{array}
\right),  \ \ f_2= \left(
\begin{array}{c}
x_0x_3 \\
x_0x_1+x_2x_3 \\
x_0^2
\end{array}
\right),  \ \ f_3= \left(
\begin{array}{c}
x_1x_3 \\
x_1^2 \\
x_0x_1-x_2x_3
\end{array}
\right)
\]
This can be obtained by transposing the matrix that appears in the
first step of a free resolution of $\beta_0^T$ over $Q$ (i.e. find
$(\Syz(\beta_0^T))^T$ over $Q$). Let $F'$ be the module generated
by $f_2,f_3$. The syzygy module $\Syz(f_2,f_3)$ over $Q$ is
generated by the columns of the matrix
\[
\left(
 \begin{array}{ccc}
  -x_1x_3 & -x_1^2 & -x_0x_1+x_2x_3 \\
  x_0x_3 & x_0x_1+x_2x_3 & x_0^2
 \end{array}
\right).
\]
Let $N$ be the extension module of $F'$ to $S$. Then $F'$ is
isomorphic to $N/x_3^2N$, and hence over $S$, $F'$ has the
presentation
\[
\gamma_0= \left(
 \begin{array}{ccccc}
 -x_1x_3 & -x_1^2 & -x_0x_1+x_2x_3 & x_3^2 & 0 \\
  x_0x_3 & x_0x_1+x_2x_3 & x_0^2  & 0  & x_3^2
 \end{array}
 \right).
\]
This corresponds to the homomorphism $s : E \rightarrow 2S$, and
hence to the injective sheaf morphism $\EE \rightarrow 2\OO$. Let
$I$ be the ideal generated by the $2 \times 2$ minors of
$\gamma_0$. Then $(I:x_3^2)$ defines the empty set in $\P^3$,
which implies by Remark~\ref{th:locally-free} that $\EE$ is a
vector bundle on $\P^3$.

By resolving $\gamma_0$, we get a minimal free resolution of the
following type for $E$:
\begin{equation} 0\rightarrow S(-4) \rightarrow 4S(-3) \rightarrow
5S(-2) \rightarrow E \rightarrow 0 \label{sq:minimal-free-res}
\end{equation}
From Sequence~(\ref{sq:minimal-free-res}) it follows that the
Chern classes of $\widetilde{E}$ are $c_1=-2$ and $c_2=2$. So the
corresponding normalized bundle is a stable rank two vector bundle
on $\P^3$ with Chern classes $(c_1,c_2)=(0,1)$.
\label{th:null-correlation}
\end{ex}
\begin{rem}
A construction almost identical to the one outlined in the previous
example can be carried out with $\Omega_{\P^n}(2)$ whenever $n$ is
even. The construction yields a rank $n$ bundle on $\P^{n+1}$.
\end{rem}

In the next example, we will discuss the stable rank two vector
bundle $\EE$ on $\P^4$ over an algebraically closed field $K$ of
characteristic two constructed by Kumar~\cite{kumar}. He proved the
existence of this bundle by constructing a rank three vector bundle
on $\P^3$ over $K$ that satisfies the conditions in
Theorem~\ref{th:kumar-corr}. Our main goal is to describe $\EE$
explicitly by using the algorithm.

\begin{ex}
Let $K$ be an algebraically closed field with characteristic two,
let $R=K[x_0,\dots,x_3]$  and let $S=K[x_0,\dots,x_4]$. Consider
the module $M$ obtained as the cokernel of the map
\[
 \alpha_0=
 \left(
  \begin{array}{cccc}
   0 & 0 & x_0x_1^2 & x_1^3 \\
   0 & 0 & x_0^3 & x_0^2x_1 \\
   x_2^2 & x_3^2 & 0 & 0 \\
   x_0 & 0 & 0 & x_3^2 \\
   0 & x_1 & x_2^2 & 0 \\
   x_1 & x_0 & x_3^2 & x_2^2
  \end{array}
 \right)
 :
 2R(-4)\oplus2R(-5) \rightarrow 3R(-2)\oplus3R(-3).
\]
Let $I_i(M)$ be the ideal of $i\times i$ minors of $\alpha_0$ (i.e.
a Fitting invariant of $M$). Then $\sqrt{I_3(M)}=(1)$ and
$I_4(M)=0$. By Fitting's Lemma, the corresponding coherent sheaf
$\widetilde{M}$ is a rank three vector bundle on $\P^3$. Let
$\phi_0$ be the homomorphism from $3R(-3)\oplus3R(-4)$ to
$3R(-2)\oplus3R(-3)$ given by
\[
\left(
 \begin{array}{cccccc}
  0 & 0 & 0 & 0 & 0 & 0 \\
  0 & 0 & 0 & 0 & 0 & 0 \\
  0 & 0 & 0 & x_1^2 & x_0^2 & x_0x_1 \\
  0 & 1 & 0 & 0 & 0 & 0 \\
  1 & 0 & 0 & 0 & 0 & 0 \\
  0 & 0 & 0 & 0 & 0 & 0
  \end{array}
\right).
\]
It is easy to check that $\phi_0$ induces a nilpotent endomorphism
$\phi$ of $M$, and hence of $\widetilde{M}$, whose third power is
zero. Therefore, $M$ corresponds to an $S$-module $F$ in
$\mathfrak{S}_3$. Let $\mathfrak{M}=\{g_i\}_{1\leq i \leq 6}$ be a
minimal generating set of $M$. Then $F$ is obtained as the
following set:
\[
 F=\{ \ a_1g_1+\cdots+a_6g_6 \ | \ a_i \in S \ \mbox{for each $i=1, \dots, 6$} \ \}
\]
The relations among $g_i$'s in $S$ are, by Proposition~\ref
{th:free_presentation}, given by the matrix
\[
\left(
%\begin{array}{cc}
\alpha_0 \arrowvert \ \phi_0-\cdot x_4
%\end{array}
\right)=
 \left(
 \begin{array}{cccccccccc}
 0 & 0 & x_0x_1^2 & x_1^3 & x_4 & 0 & 0 & 0 & 0 & 0 \\
 0 & 0 & x_0^3 & x_0^2x_1 & 0 & x_4 & 0 & 0 & 0 & 0 \\
 x_2^2 & x_3^2 & 0 & 0 & 0 & 0 & x_4 & x_1^2 & x_0^2 & x_0x_1 \\
 x_0 & 0 & 0 & x_3^2 & 0  & 1 & 0 & x_4 & 0 & 0 \\
 0 & x_1 & x_2^2 & 0 & 1 & 0 & 0 & 0 & x_4 & 0 \\
 x_1 & x_0 & x_3^2 & x_2^2 & 0 & 0 & 0 & 0 & 0 & x_4
\end{array}
\right).
\]
From the ``ones'' in this matrix, it follows that the minimal set of
generators for $F$ consists of $g_1,g_2,g_3$ and $g_6$. Eliminating
the redundant elements $g_4$ and $g_5$, we obtain a minimal free
presentation of $F$:
\[
S(-3)\oplus 5S(-4) \oplus 2S(-5)
\stackrel{\beta_0}{\longrightarrow}  3S(-2)\oplus S(-3)
\rightarrow F \rightarrow 0
\]
where
\[
\beta_0= \left(
\begin{array}{cccccccc}
0 & x_4^2 & 0 & 0 & 0 & x_1x_4 & x_0x_1^2+x_2^2x_4 & x_1^3 \\
0 & 0 & 0 & x_4^2 & x_0x_4 & 0 & x_0^3 & x_0^2x_1+x_3^2x_4 \\
x_4 & x_0^2 & x_0x_1 & x_1^2 & x_2^2 & x_3^3 & 0 & 0 \\
0 & 0 & x_4 & 0  & x_1 & x_0 & x_3^2 & x_2^2
\end{array}
\right).
\]

Next we define a homomorphism $\psi_0$ from $2R(-2)$ to
$3R(-2)\oplus3R(-3)$ by
\[
  \left(
 \begin{array}{cccccc}
 1 & 0 & 0 & 0 & 0 & 0 \\
 0 & 1 & 0 & 0 & 0 & 0
 \end{array}
 \right)^T.
\]
This represents a homomorphism $\psi$ from $2R(-2)$ to $M$. The
cokernel $C$ of $(\phi,\psi)$ has the presentation matrix
$(\phi_0,\psi_0,\alpha_0)$.
%\[
% 0 \leftarrow C \leftarrow  3R(-2)\oplus3R(-3) \stackrel{(\phi_0,\psi_0,\beta_0)}{\longleftarrow}
%3R(-3)\oplus3R(-4)\oplus2R(-2) \oplus 2R(-4)\oplus2R(-5).
%\]
Minimizing the generators and the corresponding relations, we
obtain the following presentation matrix of $C$:
\[
 \left(
 \begin{array}{ccccccc}
 x_0^2 & x_0x_1 & x_1^2 & x_2^2 & x_3^2 & 0 & 0 \\
 0 & 0 & 0 & x_1 & x_0 & x_2^2 & x_3^2
 \end{array}
 \right).
\]
Clearly $C$ is an $R$-module of finite length. From
Theorem~\ref{th:kumar-corr} it follows that there exist a rank two
vector bundle $\EE$ on $\P^4$ and an exact sequence
\begin{eqnarray}
 0 \rightarrow \EE \rightarrow 2\OO(-2) \rightarrow \widetilde{F} \rightarrow 0.
\label{sq:kumar-bundle}
\end{eqnarray}
By construction, the surjective map $2\OO(-2) \rightarrow
\widetilde{F}$ in Sequence~(\ref{sq:kumar-bundle}) is defined by
$g_1$ and $g_2$. Let $Q$ be the quotient ring $S/(x_4^3)$. Let
$P=x_0^4x_1^2+x_0^2x_1x_3^2x_4+x_3^4x_4^2
x_0x_1^5+x_2^2x_4(x_0^3+x_1^3)$. Then $\Syz(g_1,g_2)$ is generated
by the columns of the matrix $\gamma_0$ with:
\[
\gamma_0^T=\left(
\begin{array}{cc}
x_0^2x_4^2 & x_1^2x_4^2\\
x_0^4x_4 & x_0^2x_1^2x_4+x_0x_2^2x_4^2+x_1x_3^2x_4^2\\
x_0^3x_1x_4+x_0x_3^2x_4^2 & x_0x_1^3x_4+x_1x_2^2x_4^2\\
x_0^2x_1^2x_4+x_0x_2^2x_4^2+x_1x_3^2x_4^2 & x_1^4x_4\\
x_0^3x_2^2+x_0^2x_1x_3^2+x_3^4x_4 &
x_0x_1^2x_2^2+x_1^3x_3^2+x_2^4x_4\\
x_0^4x_1^2+x_0^3x_2^2x_4+x_0^2x_1x_3^2x_4 &
x_0^2x_1^4+x_2^4x_4^2\\
x_0^5x_1+x_0^3x_3^2x_4 &
x_0^3x_1^3+x_0^2x_1x_2^2x_4+x_2^2x_3^2x_4^2\\
x_0^6 & x_0^3x_1^3+x_0^2x_1x_2^2x_4+x_2^2x_3^2x_4^2\\
x_0^3x_1^3+x_0x_1^2x_3^2x_4+x_2^2x_3^2x_4^2 &
P\\
x_0^2x_1^4+x_0x_1^2x_2^2x_4+x_1^3x_3^2x_4+x_2^4x_4^2 & x_1^6\\
\end{array}
\right).
\]

%\[
%\gamma_0= \left(
%\begin{array}{cccccccccc}
%x_0^2x_4^2  & x_0^4x_4  & x_0^3x_1x_4+x_0x_3^2x_4^2 &
%x_0^2x_1^2x_4+x_0x_2^2x_4^2+x_1x_3^2x_4^2 &
%x_0^3x_2^2+x_0^2x_1x_3^2+x_3^4x_4 &
%x_0^4x_1^2+x_0^3x_2^2x_4+x_0^2x_1x_3^2x_4 & x_0^5x_1+x_0^3x_3^2x_4
%& x_0^6  &
%x_0^3x_1^3+x_0x_1^2x_3^2x_4+x_2^2x_3^2x_4^2 &  x_0^2x_1^4+x_0x_1^2x_2^2x_4+x_1^3x_3^2x_4+x_2^4x_4^2 \\
%x_1^2x_4^2 & x_0^2x_1^2x_4+x_0x_2^2x_4^2+x_1x_3^2x_4^2  &
%x_0x_1^3x_4+x_1x_2^2x_4^2 & x_1^4x_4 &
%x_0x_1^2x_2^2+x_1^3x_3^2+x_2^4x_4 & x_0^2x_1^4+x_2^4x_4^2 &
%x_0^3x_1^3+x_0^2x_1x_2^2x_4+x_2^2x_3^2x_4^2 &
%x_0^4x_1^2+x_0^3x_2^2x_4+x_0^2x_1x_3^2x_4+x_3^4x_4^2
%x_0x_1^5+x_1^3x_2^2x_4 & x_1^6
%\end{array}
%\right).
%\]
Let $N$ denote the extension module of $\Syz(g_1,g_2)$ to $S$.
Since $\Syz(g_1,g_2)$ can be identified with $N/x_4^3N$,
$\Syz(g_1,g_2)$ has, as an $S$-module, the minimal free
presentation $ \gamma= \left(
\begin{array}{cc}
\gamma_0 & \gamma_1
\end{array}
\right), $ where
\[
 \gamma_1=
 \left(
 \begin{array}{cc}
 x_4^3 & 0 \\
 0 & x_4^3
 \end{array}
 \right).
\]
This corresponds to an injective sheaf morphism $s : \EE
\rightarrow 2\OO(-2)$, whose cokernel equals $\widetilde{F}$.

Let $I$ be the ideal generated by $2 \times 2$ minors of $\gamma$.
Then the ideal quotient $(I:x_4^3)$ defines the empty subset of
$\P^4$. By Remark~\ref{th:locally-free}, $\EE$ is a rank two
vector bundle on $\P^4$. The Chern classes of $\EE$ are $c_1=-7$
and $c_2=16$. These can be computed in the same way as in
Example~\ref{th:null-correlation}.

\label{th:kumar-bundle}
\end{ex}

 As a final example, we will illustrate how to determine the triple $(\MM, \phi, \psi)$ from the pair
 $(\EE,s)$. In general, this direction is easier to carry out with the main difficulty coming from
 producing the pair $(\EE,s)$. We will discuss the Horrocks-Mumford
 bundle utilizing the ideas of Kaji to produce the sections $s$ required in the correspondence \cite{Kaji}.

\begin{ex} Let $V$ be a five-dimensional vector space with
basis $\{e_0, \dots, e_4\}$ over $K$, let $W$ be its dual and let
$\P^4=\P(V)$ be the projective space of lines in $V$. The
homogeneous coordinate ring $K[x_0, \dots, x_4] $ of $\P^4$ will be
denoted by $S$. Consider the Koszul  complex resolving $K=S/\langle
W \rangle$:
\[0\rightarrow \bigwedge^5 W \otimes S(-5)
\stackrel{\beta_4}{\rightarrow} \cdots
\stackrel{\beta_1}{\rightarrow} \bigwedge^1 W \otimes S(-1)
\stackrel{\beta_0}{\rightarrow} \bigwedge^0 W \otimes S \rightarrow
K \rightarrow 0\]

Recall that the $i^{th}$ bundle of differentials
$\Omega^i=\Omega^i_{\P^4}$ is obtained as a sheafication of the
syzygy module $\Syz_{i+1}(K)$.
%Fix the ordered bases for $\bigwedge^2 V$ and $\bigwedge^3 V$
%(and $\bigwedge^2 W$ and $\bigwedge^3 W$) respectively:
%\[
%\left\{e_0 \wedge e_1, e_0 \wedge e_2, e_1 \wedge e_2,
%e_0 \wedge e_3, e_1 \wedge e_3, e_2 \wedge e_3,
%e_0 \wedge e_4, e_1 \wedge e_4, e_2 \wedge e_4,e_3 \wedge e_4\right\}
%\]
%and
%\begin{eqnarray*}
%\left\{
%e_0 \wedge e_1 \wedge e_2, e_0 \wedge e_1 \wedge e_3, e_1 \wedge e_2 \wedge e_3,
%e_0 \wedge e_1 \wedge e_4, e_0 \wedge e_2 \wedge e_4,\right.\\
%\left.  e_1 \wedge e_2 \wedge e_4,
%e_0 \wedge e_3 \wedge e_4, e_1 \wedge e_3 \wedge e_4, e_2 \wedge e_3 \wedge e_4,
%e_2 \wedge e_3 \wedge e_4
%\right\}.
%\end{eqnarray*}
By choosing appropriate bases for $\bigwedge^2 W$ and $\bigwedge^3
W$, we may suppose that  $\mathrm{Syz}_3(K)$ is generated by the
columns of the following matrix:
\[
\beta_2=\left(
\begin{array}{cccccccccc}
{{x}}_{{2}}&
      {{x}}_{{3}}&
      0&
      0&
      {{x}}_{{4}}&
      0&
      0&
      0&
      0&
      0\\
      {-{{x}}_{1}}&
      0&
      {{x}}_{{3}}&
      0&
      0&
      {{x}}_{{4}}&
      0&
      0&
      0&
      0\\
      {{x}}_{0}&
      0&
      0&
      {{x}}_{{3}}&
      0&
      0&
      {{x}}_{{4}}&
      0&
      0&
      0\\
      0&
      {-{{x}}_{1}}&
      {-{{x}}_{{2}}}&
      0&
      0&
      0&
      0&
      {{x}}_{{4}}&
      0&
      0\\
      0&
      {{x}}_{0}&
      0&
      {-{{x}}_{{2}}}&
      0&
      0&
      0&
      0&
      {{x}}_{{4}}&
      0\\
      0&
      0&
      {{x}}_{0}&
      {{x}}_{1}&
      0&
      0&
      0&
      0&
      0&
      {{x}}_{{4}}\\
      0&
      0&
      0&
      0&
      {-{{x}}_{1}}&
      {-{{x}}_{{2}}}&
      0&
      {-{{x}}_{{3}}}&
      0&
      0\\
      0&
      0&
      0&
      0&
      {{x}}_{0}&
      0&
      {-{{x}}_{{2}}}&
      0&
      {-{{x}}_{{3}}}&
      0\\
      0&
      0&
      0&
      0&
      0&
      {{x}}_{0}&
      {{x}}_{1}&
      0&
      0&
      {-{{x}}_{{3}}}\\
      0&
      0&
      0&
      0&
      0&
      0&
      0&
      {{x}}_{0}&
      {{x}}_{1}&
      {{x}}_{{2}}\\
\end{array}
\right).
\]
The natural duality $\bigwedge^p V \otimes \bigwedge^p W \rightarrow
K$ extends to a contraction map
\[
 \bigwedge^p V \otimes \bigwedge^q W \rightarrow  \left\{
 \begin{array}{cc}
 \bigwedge^{p-q} V & \mbox{if $p \geq q$} \\
  \bigwedge^{q-p} W & \mbox{otherwise}.
 \end{array}
 \right.
\]
Using this, the linear transformation
\[
 \left(
 \begin{array}{ccccc}
      {{e}}_{{2}} \wedge {{e}}_{{3}}&
      {{e}}_{0} \wedge{{e}}_{{4}}&
      {{e}}_{1} \wedge{{e}}_{{2}}&
      {-{{e}}_{{3}} \wedge{{e}}_{{4}}}&
      {{e}}_{0}\wedge {{e}}_{1}\\
      {{e}}_{1} \wedge{{e}}_{{4}}&
      {{e}}_{1} \wedge{{e}}_{{3}}&
      {{e}}_{0}\wedge {{e}}_{{3}}&
      {{e}}_{0} \wedge{{e}}_{{2}}&
      {-{{e}}_{{2}}\wedge {{e}}_{{4}}}\\
       \end{array}
   \right)
\]
from $5 \bigwedge^5 W$ to $2 \bigwedge^2 W$ induces a sheaf
morphism, $A$, from $5 \bigwedge^5 W \otimes \OO(-1)$ to $ 2
\Omega^2(2)$. The matrix representation $A_0$ of this morphism with
respect to the fixed bases for $\bigwedge^2 W$ and $\bigwedge^3 W$
is
\[
{\tiny
 \left(
 \begin{array}{cccccccccccccccccccc}
0&
      0&
      0&
      0&
      {-1}&
      0&
      0&
      0&
      0&
      0&
      0&
      0&
      {-1}&
      0&
      0&
      0&
      0&
      0&
      0&
      0\\
      0&
      0&
      0&
      1&
      0&
      0&
      0&
      0&
      0&
      0&
      0&
      0&
      0&
      0&
      0&
      1&
      0&
      0&
      0&
      0\\
      0&
      0&
      0&
      0&
      0&
      0&
      0&
      {-1}&
      0&
      0&
      0&
      0&
      0&
      0&
      0&
      0&
      {-1}&
      0&
      0&
      0\\
      1&
      0&
      0&
      0&
      0&
      0&
      0&
      0&
      0&
      0&
      0&
      0&
      0&
      0&
      0&
      0&
      0&
      0&
      1&
      0\\
      0&
      0&
      0&
      0&
      0&
      0&
      0&
      0&
      0&
      {-1}&
      0&
      {-1}&
      0&
      0&
      0&
      0&
      0&
      0&
      0&
      0\\
 \end{array}
 \right)^T
 }.
\]
Let $\beta=
 \left(\begin{array}{cc}
 \beta_2 & 0 \\
 0 & \beta_2
 \end{array}
 \right)$.
One can show that the ideal generated  by the maximal 
minors of the composite of  $\beta$ and $A_0$ defines the empty
set, and thus $A$ is injective as a bundle map. Let $B_0=A_0^T \cdot
\left(\begin{array}{cc} 0 & I_{10} \\ -I_{10} & 0 \end{array}
\right)$ where $I_{10}$ is the $10 \times 10$ identity matrix. The
matrix $B_0$ gives rise to a sheaf morphism  $B$ from $2\Omega^2(2)$
to $5\bigwedge^0 W \otimes \OO$. This sheaf morphism is surjective
as a bundle map (since $A$ is injective). $A$ and $B$ can be thought
of as the differentials of the following complex:
\[
5 \bigwedge^5 W \otimes \OO(-1) \stackrel{A}{\rightarrow}
2 \Omega^2 (2) \stackrel{B}{\rightarrow} 5 \bigwedge^0 W \otimes \OO.
\]
Since $A$ is an injective bundle map and $B$ is a surjective bundle
map the homology, $\EE=\ker{B}/\im{A}$, is a rank two vector bundle
on $\P^4$. This vector bundle is known as the Horrocks-Mumford
bundle, is indecomposable and has Chern classes $c_1=-1$ and
$c_2=4$.

% \left(
%\begin{array}{ccccc}
%     {-{{e}}_{1} \wedge{{e}}_{{4}}}&
%     {{e}}_{0} \wedge{{e}}_{{2}}&
%     {{e}}_{1} \wedge{{e}}_{{3}}&
%     {{e}}_{{2}} \wedge{{e}}_{{4}}&
%     {-{{e}}_{0} \wedge{{e}}_{{3}}}\\
%     {{e}}_{{2}} \wedge{{e}}_{{3}}&
%     {{e}}_{{3}} \wedge{{e}}_{{4}}&
%     {-{{e}}_{0}\wedge {{e}}_{{4}}}&
%     {{e}}_{0} \wedge{{e}}_{1}&
%     {{e}}_{1} \wedge{{e}}_{{2}}\\
% \end{array}
% \right)^T : 2 \bigwedge^2 W \rightarrow 5 \bigwedge^0 W
%\]
Consider the following $20 \times 1$ matrices $v_1$ and $v_2$
(discovered by Kaji \cite{Kaji})
\[v_1= {\footnotesize
 \left(
 %\begin{array}{cccccccccccccccccccc}
      0, 
      B_2.
      0, 
      0, 
      0, 
      B_6, 
      0, 
      0, 
      0, 
      0, 
      B_{11}, 
      B_{12}, 
      0, 
      0, 
      B_{15},  
      0, 
      0, 
      B_{18}, 
      0, 
      0\\
% \end{array}
 \right)^T
 }
\]
\[v_2= {\footnotesize
 \left(
 \begin{array}{cccccccccccccccccccc}
      0, 
      C_2, 
      C_3, 
      0, 
      0, 
      C_6, 
      0, 
      0, 
      0, 
      0, 
      C_{11}, 
      C_{12}, 
      0, 
      0, 
      0, 
      0, 
      0, 
      C_{18}, 
      0, 
      0\\
 \end{array}
 \right)^T
 }
\] where

\[ \begin{array}{lll}
B_2 = -{{x}}_{0}^{5} {{x}}_{1}-{{x}}_{0} {{x}}_{1}^{2} {{x}}_{{2}}
{{x}}_{{3}} {{x}}_{{4}}-{{x}}_{0}^{3} {{x}}_{{3}} {{x}}_{{4}}^{2} &
C_2 =
 -{{x}}_{0}^{5} {{x}}_{{3}}^{2}-{{x}}_{0}^{3} {{x}}_{{2}}^{2} {{x}}_{{3}}
 {{x}}_{{4}}-{{x}}_{0} {{x}}_{1} {{x}}_{{2}} {{x}}_{{3}}^{3} {{x}}_{{4}}\\
B_6 =
 -{{x}}_{0}^{3} {{x}}_{1}^{2} {{x}}_{{2}}-{{x}}_{0}^{5} {{x}}_{{4}}  & C_3
 =
 {-{{x}}_{0}^{7}}\\
B_{11} =
 {{x}}_{0}^{4} {{x}}_{1}^{2}+{{x}}_{1}^{3} {{x}}_{{2}} {{x}}_{{3}}
 {{x}}_{{4}}+{{x}}_{0}^{2} {{x}}_{1} {{x}}_{{3}} {{x}}_{{4}}^{2} & C_6
 =
-{{x}}_{0}^{5} {{x}}_{{2}}^{2}-{{x}}_{0}^{3} {{x}}_{1} {{x}}_{{2}} {{x}}_{{3}}^{2}\\
B_{12} =
 -{{x}}_{1}^{3} {{x}}_{{2}}^{2} {{x}}_{{4}}-{{x}}_{0}^{2} {{x}}_{1} {{x}}_{{2}} {{x}}_{{4}}^{2} & C_{11}
 =
{{x}}_{0}^{6} {{x}}_{{2}}+{{x}}_{0}^{4} {{x}}_{1} {{x}}_{{3}}^{2}+{{x}}_{0}^{2}
{{x}}_{1} {{x}}_{{2}}^{2} {{x}}_{{3}} {{x}}_{{4}}\\
& \hskip 32pt +{{x}}_{1}^{2} {{x}}_{{2}} {{x}}_{{3}}^{3} {{x}}_{{4}}\\
B_{15} =
 {{x}}_{0}^{6} &  C_{12} =
-{{x}}_{0}^{2} {{x}}_{1} {{x}}_{{2}}^{3} {{x}}_{{4}}-{{x}}_{1}^{2} {{x}}_{{2}}^{2} {{x}}_{{3}}^{2} {{x}}_{{4}}\\
B_{18} =
 {{x}}_{0}^{2} {{x}}_{1}^{2} {{x}}_{{2}} {{x}}_{{4}}+{{x}}_{0}^{4}
 {{x}}_{{4}}^{2} & C_{18} =
 {{x}}_{0}^{6} {{x}}_{{3}}+{{x}}_{0}^{4} {{x}}_{{2}}^{2} {{x}}_{{4}}+{{x}}_{0}^{2}
 {{x}}_{1} {{x}}_{{2}} {{x}}_{{3}}^{2} {{x}}_{{4}}.\\
\end{array}\]

The matrix $v_1$ represents a global section $s_1$ of
$2\Omega^2(9)$; %\simeq \Hom(\OO(-7), 2\Omega^2(2))$;
while $v_2$ represents a global section $s_2$ of
$2\Omega^2(10)$. % \simeq \Hom(\OO(-8), 2\Omega^2(2))$.
Both $v_1$ and $v_2$ can be  written as $S$-linear combinations of
the columns of $\Syz(\beta \circ B_0)$, thus $s_1$ and $s_2$
correspond to global sections $\widetilde{s}_1$ and
$\widetilde{s}_2$ of $\EE(7)$ and $\EE(8)$ respectively. Both
$\widetilde{s}_1$ and $\widetilde{s}_2$ are nonzero and together
generate $\EE$ on $D_+(x_0)$. Indeed, if $I$ is the ideal generated
by the maximal minors of the matrix $\left(
%\begin{array}{ccc}
v_1, v_2, A_0
% \end{array}
\right)$ then the saturation of $I$ with respect to $x_0$ determines
the locus of points, not on $H$, where $s_1$ and $s_2$ do not
generate $\EE$ ($H$ is the hyperplane defined by $x_0=0$). An easy
computation establishes that $V(I:(x_0)^{\infty})=V((1))=\emptyset$.

The global sections $\widetilde{s}_1$ and  $\widetilde{s}_2$ can be
identified with a sheaf morphism $s=(s_1,s_2)$ from $\OO(-8) \oplus
\OO(-7)$ to $\EE$. Recall that $\EE^{\vee}$ is isomorphic to
$\EE(c_1)$ (since $\EE$ is a rank 2 reflexive sheaf). Taking the
transpose of $s$ we obtain the following short exact sequence:
\[
 0 \rightarrow \EE(-1) \stackrel{s^{\vee}}{\rightarrow}
 \OO(7) \oplus \OO(8)
 \rightarrow \FF \rightarrow 0.
\]
Since $\widetilde{s}_1 \wedge \widetilde{s}_2 \in \H^0(\P^4, \EE(7)
\wedge \EE(8)) \simeq  \H^0(\P^4, \OO(14))$ and since $s_1,s_2$
generate $\EE$ away from $H$, the sheaf $\FF$ can be considered as a
coherent sheaf on the $14$th infinitesimal neighborhood $H_{14}$ of
$H$. Let $\pi$ be the finite morphism from $H_{14}$ to $H$ induced
by the projection $\P^4 \setminus P \rightarrow H$ from a point $P$
off $H$. Then the direct image sheaf of $\FF$ by $\pi$ is a rank
fourteen vector bundle on $H\simeq \P^3$. We denote this bundle by
$\MM$.

Let $R$ be  the quotient ring %$S/\left\langle x_0^{14} \right\rangle$ and
$S/(x_0)$ and let $F$ be the graded $T$-module $\H^0_* \FF$. Then
the graded $R$-module $M=\H^0_* \MM$ is the graded $R$-module
${}_RF$ obtained from $F$ by restriction of scalars. It is
straightforward to determine that $F$ has a minimal free
presentation of the following form:
\[
% 0 \leftarrow F \leftarrow S(7) \oplus S(8) \oplus 5S(1)
% \stackrel{A}{\leftarrow} 15 S.
 15S \stackrel{P}{\longrightarrow} S(8) \oplus S(7) \oplus 5S(1)  \rightarrow F \rightarrow 0,
\]
 Let
$\mathfrak{F}=\{f_i\}_{1 \leq i \leq 7}$ be the minimal generating
set of $F$. Then it follows from Proposition 2.2 that $\mathfrak{M}=
\left\{ \ x_0^if_j  \ | \ 0\leq i \leq 13, 1 \leq j \leq 7 \
\right\}
%_{{0 \atop 1}{\leq \atop \leq }{i \atop j}{\leq \atop \leq}{13 \atop 7}}
$ is a set of generators for $M$. The relations among these
generators of $M$ can be derived from the presentation matrix $P$ of
$F$. Let $P[:,k]$ be the $k^{th}$ column of $P$ and let $Q$ be the
presentation matrix of $M$ with respect to $\mathfrak{M}$. For each
$1\leq k \leq 15$, we have a syzygy of the form
\[
\sum_{i=1}^{7}P[i,k]f_i=0.
\]
Then, since
\[
P[i,k]=\sum_{t=0}^{13}Q[7t+i,k]x_0^t
\]
we can obtain the entries of $Q[:,k]$ from the entries of $P[:,k]$.

Choosing appropriate bases for $F_0$ and $F_1$, one can explicitly
write $P$. For example, the first column of $P$ is
\[
 P[:,1]=
 \left(
\begin{array}{ccccccc}
P[1,1] & P[2,1] &
 x_3  &
      0 &
      0 &
      0 &
      0
\end{array}
\right)^T,
\]
where
\begin{eqnarray*}
 P[1,1] & = & {{x}}_{0}^{6} {{x}}_{{2}}^{2}-{{x}}_{1}^{3} {{x}}_{{2}}^{4} {{x}}_{{3}}+2 {{x}}_{0}^{2} {{x}}_{1} {{x}}_{{2}}^{3} {{x}}_{{3}} {{x}}_{{4}}+{{x}}_{0} {{x}}_{1}^{4} {{x}}_{{3}}^{2} {{x}}_{{4}}-3 {{x}}_{1}^{2} {{x}}_{{2}}^{2} {{x}}_{{3}}^{3} {{x}}_{{4}} \\
 &  &
 -{{x}}_{0}^{2} {{x}}_{{2}} {{x}}_{{3}}^{3} {{x}}_{{4}}^{2}
 +{{x}}_{1} {{x}}_{{3}}^{5} {{x}}_{{4}}^{2}-{{x}}_{0} {{x}}_{1}^{2} {{x}}_{{2}} {{x}}_{{3}} {{x}}_{{4}}^{3}-{{x}}_{{2}}^{3} {{x}}_{{3}}^{2} {{x}}_{{4}}^{3}+{{x}}_{0}^{3} {{x}}_{{3}} {{x}}_{{4}}^{4}, \\
P[2,1] & = & {{x}}_{0}^{4} {{x}}_{1}^{2} {{x}}_{{2}}-{{x}}_{0}^{3} {{x}}_{{2}}^{3} {{x}}_{{3}}+{{x}}_{0} {{x}}_{1} {{x}}_{{2}}^{2} {{x}}_{{3}}^{3}+{{x}}_{0}^{6} {{x}}_{{4}}+{{x}}_{1}^{3} {{x}}_{{2}}^{2} {{x}}_{{3}} {{x}}_{{4}}-{{x}}_{0} {{x}}_{{3}}^{5} {{x}}_{{4}} \\
&  & +2 {{x}}_{0}^{2} {{x}}_{1} {{x}}_{{2}} {{x}}_{{3}}
{{x}}_{{4}}^{2}+{{x}}_{1}^{2} {{x}}_{{3}}^{3}
{{x}}_{{4}}^{2}-{{x}}_{{2}} {{x}}_{{3}}^{2} {{x}}_{{4}}^{4}.
\end{eqnarray*}
\[
{\rm We\ have}\ \ P[1,1]=Q[1,1]+Q[8,1] x_0 +Q[15,1] x_0^2+Q[22,1]
x_0^3 +Q[36,1] x_0^6,
\]
where
\begin{eqnarray*}
Q[1,1] & = & -{{x}}_{1}^{3} {{x}}_{{2}}^{4} {{x}}_{{3}}-3 {{x}}_{1}^{2} {{x}}_{{2}}^{2} {{x}}_{{3}}^{3} {{x}}_{{4}}+{{x}}_{1} {{x}}_{{3}}^{5} {{x}}_{{4}}^{2}-{{x}}_{{2}}^{3} {{x}}_{{3}}^{2} {{x}}_{{4}}^{3}\\
Q[8,1] & = & {{x}}_{1}^{4} {{x}}_{{3}}^{2} {{x}}_{{4}}-{{x}}_{1}^{2} {{x}}_{{2}} {{x}}_{{3}} {{x}}_{{4}}^{3}\\
Q[15,1] & = & 2 {{x}}_{1} {{x}}_{{2}}^{3} {{x}}_{{3}} {{x}}_{{4}}-{{x}}_{{2}} {{x}}_{{3}}^{3} {{x}}_{{4}}^{2}\\
Q[22,1] & = &  {{x}}_{{3}} {{x}}_{{4}}^{4} \\
Q[36,1] & = & {{x}}_{{2}}^{2}
\end{eqnarray*}
Likewise, $$P[2,1]=Q[2,1]+Q[9,1] x_0 +Q[16,1] x_0^2+Q[23,1]
x_0^3+Q[30,1] x_0^4+Q[37,1]x_0^6$$ where
\begin{eqnarray*}
Q[2,1] & = & {{x}}_{1}^{3} {{x}}_{{2}}^{2} {{x}}_{{3}} {{x}}_{{4}}+{{x}}_{1}^{2} {{x}}_{{3}}^{3} {{x}}_{{4}}^{2}-{{x}}_{{2}} {{x}}_{{3}}^{2} {{x}}_{{4}}^{4}\\
Q[9,1] &=& {{x}}_{1} {{x}}_{{2}}^{2} {{x}}_{{3}}^{3}-{{x}}_{{3}}^{5} {{x}}_{{4}}\\
Q[16,1] &=&  2 {{x}}_{1} {{x}}_{{2}} {{x}}_{{3}} {{x}}_{{4}}^{2}\\
Q[23,1] &=& {-{{x}}_{{2}}^{3} {{x}}_{{3}}}\\
Q[30,1] &=& {{x}}_{1}^{2} {{x}}_{{2}}\\
Q[37,1] &=& {{x}}_{{4}}.
\end{eqnarray*}
Finally, $Q[3,1]=x_3$ is the remaining nonzero entry in $Q[:,1]$
(since $P[i,1]=0$ for $4 \leq i \leq  7$).

Working our way through the other columns of $P$, the entire matrix
$Q$ can be obtained (and has $98=14 \cdot 7$ rows and $15$ columns).
Upon obtaining $Q$, one finds that
$Q[12,6],Q[10,8],Q[11,10],Q[13,14],Q[14,15], Q[51,12]$ and  $Q[57,13]$ are
the only entries of $Q$ which are constant and nonzero. Furthermore,
each of $\{x_0^if_j|i \geq 1,j\geq 3\}$, $\{x_0^if_1|i \geq 8\}$ and
$\{x_0^if_2|i \geq 7\}$, can be written as $R$-linear combinations
of $${\mathbf G}=\{f_1,f_2, \dots, f_7\}\cup \{ x_0^if_1|1\leq i\leq
7\}\cup\{x_0^if_2|1\leq i\leq 6\}.$$ These linear combinations give
rise to the standard nilpotent endomorphism of $M$. Let $g_j$ denote
the $j^{th}$ entry of ${\mathbf G}$ for $1 \leq j \leq 20$ and let
\[
 M_0 \stackrel{\mathbf{G}}{\longrightarrow} M \rightarrow 0
\]
be the map associated to the minimal set of generators of $M$. Each
$x_0 g_i$ can be written as an $R$-linear combination of $g_1,
\dots, g_{20}$:
\[
 x_0g_i = \sum_{j=1}^{20} a_{ij} g_j.
\]
The matrix  $(a_{ij})_{1 \leq i,j  \leq 20}$ is the standard lifting
of the standard nilpotent endomorphism $\phi$ of $M$ (see Remark
$2.3$). By construction, the first two generators $g_1$ and $g_2$ of
$M$ form a homomorphism $\psi$ from $R(7) \oplus R(8)$ to $M$ such
that the cokernel of $(\phi[-1], \psi) : M(-1) \oplus R(7) \oplus
R(8) \rightarrow M$ is a finite-length $R$-module.
\end{ex}

It is interesting to note that the rank fourteen vector bundle $\MM$
can be written as the direct sum of nine line bundles and an
indecomposable rank five vector bundle.

%---reference--

\end{document}